\documentclass[a4paper,12pt]{article}
\usepackage{latexsym}
\usepackage{amsfonts}
\def \R {{\rm I\!R}}

\def \D {{\rm I\!D}}
\def \L {{\rm I\!L}}

\def \h {{\dot{h}}}
\def \e {{\bf \rm e}}
\def \u {{\bf \rm u}}

\def \tH {{\tilde{H}}}
\def \H {{\rm I\! H}}
\def \l {{\cal L}}
\def \E {{\cal E}}
\def \ric {{\rm ric_{r(s)}}}
\def \ro {\gamma}
\def \D {{\rm I\!D}}
\def \tD {{ \tilde{D}}}
\def \deltam {{ \tilde{\delta}}}

\title{On the quasi-invariance of the Wiener measure on path spaces and the 
anticipative integrals over a Riemannian manifold} 

\author{ Adnan Aboulala\^a\footnotemark[1]}
\date{}

\begin{document}
\vskip2cm
\newtheorem{thm}{Theorem}[section]
\newtheorem{defi}[thm]{Definition}
\newtheorem{prop}[thm]{Proposition}
\newtheorem{Rq}[thm]{Remark}
\newtheorem{lm}{Lemma}[section]
\newtheorem{co}[thm]{Corollary}
\maketitle
\footnotetext[1]{ 
  This paper was written in 1996, when the author was at the Laboratoire de Probabilit\'es, Universit\'e Paris VI during the period 1995-1997. Some more recent references have been added in this version.  E-mail : adnan.aboulalaa@polytechnique.org }
\begin{abstract}

\noindent
Some parts of stochastic analysis on curved spaces are revisted. A concise proof of the quasi-invariance of the Wiener measure on the path spaces over a Riemannian manifold is presented. The shifts are allowed to be in the Cameron-Martin space and random. The second part of the paper presents some remarks on the anticipative integrals on Riemannian manifolds.
\\ \\


{\it Mathematics Subject Classication} (2000): 58J65, 60H07, 28C20.
\end{abstract}

%%%%%%%%%%%%%%%%%%%%%%%%%%%%%%%%%%%%%%%%%%%%%%%%%%%%%%%%%%%%%%%%%%%%%%%%

\section{Introduction}
Let $d>0$, $W=(C_{0}([0,1], \R^{d}),\H, \mu)$ denotes the standard Wiener
 space: 
$\H=\{ h\in C_{0}([0,1], \R^{d}): \; \int_{0}^{1}|\h_{s}|^{2}ds <\infty \}$ 
and $\mu$ is 
the Wiener measure on $C_{0}([0,1], \R^{d})$ the space of $\R^{d}$-valued 
continuous functions on $[0,1]$ starting from $0$. The problem of the 
quasi-invariance
of the Wiener measure under translations $ T: w\mapsto w+k(w)$ with  
$k(w)\in \H $ a.s. 
was considered by Cameron and Martin \cite{CM} who have first dealt with the case 
where $k(w)$ is a non-random element of $\H$. In this case $T$ leaves the
Wiener measure $\mu$ quasi-invariant i.e., the measure $T_{*}\mu$ induced by 
$T$ and $\mu$ are mutually absolutely continuous. Later the case where $k$ is 
random was considered \cite{CM1} and has been 
generalized by many authors. 
The analogue of the Cameron-Martin theorem in the case 
of Brownian motion on a Riemannian manifold was first developed
 by Driver
\cite{D}. Let $M$ be a compact Riemannian manifold endowed with a connection
$\nabla $ compatible with the metric, $m_{0}\in M$ a fixed 
point and $P_{m_{0}}(M):=\{p \in C([0,1], M): p_{0}=m_{0}\}$ the path space 
over $M$. Let $\nu$ be the Wiener measure on $P_{m_{0}}(M)$ i.e. the law of
the Brownian motion $(p_{s}(w))$ on $M$ starting from $m_{0}$. A path space analogue of 
the flow $(t,w)\mapsto w+th$ which leaves the Wiener measure $\nu$ 
quasi-invariant was constructed in \cite{D}. It was shown in that paper 
that among other possible flows, {\it a} flow of semimartingales 
$\sigma_{.}(t)$ which has the quasi-invariance
property is the solution, in an appropriate sense, to
\[ \frac{d\sigma_{.}(t,w)}{dt}=t^{\sigma(t)}_{.\leftarrow 0}h_{.},\; \;
\sigma_{.}(0,w)=p_{.}(w) \]
where $t^{\sigma(t)}_{s\leftarrow 0}$ is the It\^o stochastic
parallel transport (with respect to $\nabla$) and $h\in \H$. (Note that 
$t\in \R$ is the parameter of the flow while $s$ and  ``.'' denote the time.) 
Another convenient form of this equation is
\begin{equation}
\label{eq1.1}
 \frac{d\sigma_{.}(t,w)}{dt} = H_{.}(\sigma(t,w))h_{.},
\end{equation}
where $s\mapsto H_{s}(\sigma(t))$ is the horizontal lift of the
semimartingale $\sigma(t)$. The main result of \cite{D} is that if
the function $h$ is $C^{1}$ then the 
equation (\ref{eq1.1}) has a unique solution $\sigma(t)$ in the set
of ``brownian semimartingales'', and if the torsion of
the connection is ``skew-symetric'' (see $\S$ 4 below) then, for each $t$,
the transformation $p\mapsto \sigma(t)(p)$ in the path
space $P_{m_{0}}(M)$ leaves the Wiener measure quasi-invariant. Later
 E. Hsu \cite{Hsu1} removed the restriction $C^{1}$ on $h$ to
allow it to be in the Cameron-Martin space $\H$. See also 
Enchev-Stroock \cite{ES} for another approach.

The purpose of the first part of this paper (Sections 2-4) is to prove 
the existence and uniqueness
of solutions to (\ref{eq1.1}) in the case where $h$ is random adapted
and such that  there is exists a constant 
$C$ with $\int_{0}^{1}|\h_{s}|^{2}ds \leq C$ a.s. While slight modifications
to Hsu's proof could give the same result, we feel that the proof given here
--which is inspired by that of Driver[\cite{D}, Section 7]-- is 
elementary and more direct.

Notice that in the case where $h$ is random, the transformation
$t\mapsto \sigma(t)$ will not have the flow property just like the 
flat case where the transformation $(t,w)\mapsto w+th(w)$ has not
the flow property  in general. At this point one can also
try the generalization to the case where $h$ is not adapted. This 
does not seem plausible; see remark 4.1 below. 

The ``intrinsic'' Cameron-Martin theorem for Brownian motion on Riemannian manifolds 
had renewed the interest in stochastic analysis on path and loop spaces in the mid-1990s.
In fact, let us recall that the quasi-invariance property plays
a major role in the analysis of Wiener
functionals -- the flat case -- based on the stochastic calculus of
variations; for it allows to define a closable gradient (via an integration
by parts formula) which can
be applied to a wide class of Wiener functionals; this is not possible 
with the usual differential calculus without probability theory (see e.g., \cite{Ma}, \cite{U}). Hence, 
as we have a Cameron-martin theorem on the path
space $P_{m_{0}}(M)$, one would be able to develop the same methods on 
$P_{m_{0}}(M)$. 
However, it should be noted that such a formalism can be constructed 
without using this quasi-invariance theorem, see Fang-Malliavin 
\cite{FM}. 

As in the flat case, one can seek the anticipative calculus based on 
the stochastic calculus of variations over the path space. The first step
is the study of the anticipative integrals. One of the difficulties is that a
chaotic development machinery adapted to this calculus is not
available except in the case of Lie groups (see e.g., \cite{PU}, \cite{Gr},
\cite{Us}). 
The second part of this paper contains some remarks about these anticipative
integrals. Using the parallel transport, 
a gradient $D^{M}$ on cylindrical functions on the path space $P(M)$ can be 
defined (\cite{D}, \cite{Le}, \cite{FM}).
This gradient is related to the transformation given by (\ref{eq1.1}). The
corresponding integration by part formula gives rise to a curvature term.
In order to get rid of this term, another gradient $\tD$ was 
introduced in \cite{FM}.
Thus, following Gaveau-Trauber \cite{GT}, we define two 
anticipative integrals $\delta^{M}$ and $\deltam$ which
correspond to $D^{M}$ and $\tD$ respectively. In \cite{F}, Fang proved
that the domain of $\deltam$ contains the Sobolev spaces  
$\L^{1,p}(M), p>2$ (see
Section 5). In this part we prove that the domain of $\deltam$ contains
the  spaces  $\L^{1,2}_{C}$, and
 we give an `explicit' expression of
$\deltam$ and the corresponding It\^o formula. Let us remark that the 
expression of $\delta^{M}$ is not as convenient as that of $\deltam$. We end 
this section by some remarks on the difficulties to find satisfactory $L^{p}$-estimates of
 these stochastic integrals, which are mainly due to the rotational derivative
which appear in their expressions.
\\
Some additional references: For a general account on the subject, see \cite{D1}. For another presentation and other results see Elworthy and Li \cite{EL}. Further extensions of the quasi-invariance property on path spaces may be found in Hsu \cite{Hsu2}, Hsu and Ouyang \cite{HO}, Bell \cite{Bell} and Zhang and Kannan \cite{ZK}.

%%%%%%%%%%%%%%%%%%%%%%%%%%%%%%%%%%%%%%%%%%%%%%%%%%%%%%%%%%%%%%%%%%%%%%%%%%%%%%%%%%%

\section{Preliminaries}  
\subsection{ Some geometric notations}
Let $M$ be a compact connected Riemannian manifold of dimension $d$. We denote
by $g$ the tensor metric and by ${\cal O}(M)$ the bundle of orthonormal frames 
i.e. the set of $r=(m, r_{m})$ where $m\in M$ and $r_{m}$ is an orthonormal
basis of $T_{m}M$. Throughout, we fix an element $r_{0}=(m_{0}, r_{m_{0}})$,
 $T_{m_{0}}M$ will be identified to $\R^{d}$ and for each $r=(m, r_{m})$, 
$r_{m}$ will be identified with an isometrie of $\R^{d}$ onto $T_{m}M$.

\noindent
In ${\cal O}(M)$ the canonical $\R^{d}$-valued form $\theta$ is defined
by $\theta_{u}(\xi)=r^{-1}d\pi(\xi)$, where $\pi: {\cal O}(M)\longrightarrow M$
 is the canonical projection.
We suppose that $M$ is endowed with an affine connection $\nabla$ compatible
with $g$. This connection determines an $so(d)$-valued one form $\omega$ on 
${\cal O}(M)$ by $\omega_{r}(\dot{\gamma}(0))=\gamma^{-1}(0)(\nabla
(\gamma(s))/ds(s=0))$ for every smooth curve  $\gamma$ on $O(M)$. 
 Equivalently, given $r\in {\cal O}(M)$, 
the connection $\nabla$ determines the 
horizontal subspace $H_{r}{\cal O}(M)$ of $T_{r}{\cal O}(M)$ (cf., e.g., 
 \cite{BC},
\cite{CDD}) by $\xi\in H_{r}({\cal O}(M))$ iff $\omega_{r}(\xi)=0$ 
(dim $(H_{r}({\cal O}(M))
=d$. For $v\in T_{m}M$, a vector $\tilde{v}\in T_{r}{\cal O}(M)$ is a 
horizontal lift of $v$ iif $\pi(r)=m, \tilde{v}\in H_{r}{\cal O}(M)$ and
 $d\pi(r).\tilde{v}=v$.   

Troughout, $(\e_{i}), i=1,...,d$ will denote the canonical basis of
 $\R^{d}$. To each $i\in \{1,...,d\}$ we associate a vector field $L_{i}$
on $O(M)$ as follows : for $r\in {\cal O}(M)$, $L_{i}(r)$ is the horizontal 
lift of $r\e_{i}$ ( $\omega_{r}(L_{i}(r))=0$ and $\theta_{r}(L_{i}(r))=
\e_{i}$).

Let $R, T$ be the curvature and the torsion tensors  and $\Omega$, 
$\Theta$ be the curvature ($so(d)$-valued) 2-form and the 
torsion ($\R^{d}$-valued) form respectively defined
on ${\cal O}(M)$.
 Then we have the following structural equations
\begin{equation}
\label{struct}
d\theta = -\omega\wedge\theta +\Theta,\;  d\omega=-\omega\wedge \omega +\Omega.
\end{equation} 
 Also, for $v_{1}, v_{2}\in \R^{d}$ we denote by $\Omega_{r}(v_{1},
v_{2})$ the matrix $\Omega_{r}(\tilde{r v_{1}},\tilde{r v_{2}})$ and we set
\[ {\rm ric}_{r}(v)=-\sum_{i=1}^{d}\Omega_{r}(\e_{i}, v)\e_{i},
 \; \mbox{for} \;  v\in \R^{d}\]
(the covariant representation of the Ricci tensor).

\subsection{An imbedding procedure}
In this paragraph we recall some results of [{\cite{D}, section 2] 
concerning the imbedding
of the manifold $(M,g,\nabla)$ into an open neighborhood $Y\subset \R^{N}$
 for some $N$, with a convenient extension of the covariant derivative 
$\nabla$. This extension is needed to garantee property $(i)$ of the
Proposition 2.1 below and will be used in the next section. So
let $M$ be imbedded in $\R^{N}$ for some $N$. Then :

There exists an open neighborhood $Y_{0}\subset \R^{N}$ of $M$ endowed
with a Riemannian $\bar{g}$ and a  covariant derivative 
$\bar{\nabla}$ compatible with $\bar{g}$ and a map $p: Y_{0}\rightarrow
 M\subset \R^{N}$
such that

(i) $p_{|M}=id_{M}$ and if $i: M\longrightarrow Y_{0}$ is the inclusion map 
then the 
map $j:{\cal O}(M)\longrightarrow  \R^{N}\times L(\R^{d},\R^{N})$
 $r\mapsto (i\circ \pi(r),i'(\pi(r))\circ r))$ is an imbedding
of ${\cal O}(M)$. 

(ii) If $\Gamma $ is the $N\times N$-matrix one form on $Y$
such that $\bar{\nabla}_{X} Z=d_{X}Z+\Gamma(X).Z$ for each vector
fields $X, Z$ on $Y$, then for all $y\in Y, v_{y}\in T_{y}Y$:
\begin{equation}
\label{emb}
p'(y)\Gamma(v_{y})=p''(v_{y})+\Gamma_{y}(p'(y).v_{y}).p'(y).
\end{equation}
The main consequences of this last equation is the property (i)
of proposition 2.1 below and a stochastic equation 
of the horizontal lift ( see equation (\ref{lift}) below).

Now, let $(x_{s})$ be an $M$-valued semimartingale. A semimartingale
$X_{s}$ with values in ${\cal O}(M)$ is a horizontal lift of $(x_{s})$
(with respect to the connection $\nabla$ ) if $\pi\circ X_{s}=
x_{s}$ a.s. and $\int_{[0,X_{s}]}\omega \circ dX_{\tau}$ 
vanishes a.s.(the last integral is the stochastic line integral,
see Ikeda-Watanabe \cite{IW} and the references therein).

For example the horizontal lift $r(s)$ of the Brownian motion $p_{s}$ on 
$M$ is defined by the following SDE on $O(M)$:
\[ dr(s)=L_{i}(r(s))\circ dw_{s}^{i}, \; r(0)=r_{0}. \]
This is, in fact, a method for the construction of the Brownian motion
on $M$ by setting $p_{s}=\pi(r(s))$. See \cite{IW}.

We have the following existence and uniqueness result: given 
an $M$-valued semimartingale, there is a unique horizontal
lift $X_{s}$ of $(x_{s})$ such that $X_{0}=r_{0}$. Furthermore,
 by using the above imbedding $X_{s}$ is the unique solution
to the Stratonovich SDE
\begin{equation}
\label{lift}
dX_{s}=-\Gamma_{x_{s}}(\circ dx_{s}). X_{s},\; \; X_{0}=r_{0}.
\end{equation}

Given $(y_{s})$ a $Y_{0}$-valued semimartingale, we will denote
by $(\tilde{H}_{s}(y)), s\in [0,1]$ the solution $X_{s}$ to the 
equation $dX_{s}=-\Gamma_{y_{s}}(\circ dy_{s}). X_{s}$, and 
if $y_{s}\in M$ a.s.  we denote this solution by $(H_{s}(x))$ .

In the sequel, $Y\subset Y_{0}$ will denote a compact neighborhood of $M$. So
$\Gamma$ is bounded with its derivatives on $Y$. 
\begin{prop}
Let $x_{s}, s\in [0,1]$ be a continuous $Y$-valued semimartingale such that 
$x_{0}\in M$ Then\\
(i) $p'(x_{s})\tH_{s}(x)=H_{s}(p\circ x)$.\\
(ii)if $x_{s}$ is of the form  $dx_{s}=O_{s}dw_{s}+A_{s}ds$ then 
$X_{s}=H_{s}(x)$ is of the form
\begin{equation}
dX_{s}= G_{s}(dw_{s}).X_{s}+F_{s}X_{s}ds ,
\end{equation} 
where $G_{s}$ and $F_{s}$ are such that for some constant $K$ we have:
\[
\left\{\begin{array}{c}
\sup_{s\in [0,1]}\leq K\sup_{s\in [0,1]}|O_{s}| \; \; a.s.\\
\displaystyle \int_{0}^{1}|F_{s}|^{2}ds \leq K[\int_{0}^{1}|A_{s}|^{2}ds+
\sup_{s\in[0,1]}|O_{s}|^{2}+\sup_{s\in[0,1]}|O_{s}|^{4}].
\end{array}\right. 
\]
\end{prop}

\noindent
{\bf Proof.} (i) is the lemma 7.1 of \cite{D}.

\noindent
The stochastic equation of the horizontal lift $X_{s}$ writen in the It\^o 
form is 
\begin{eqnarray}
 dX_{s}&=&\Gamma_{x_{s}}(O_{s}dw_{s}).X_{s}-\{\Gamma_{x_{s}}( A_{s})
\nonumber\\
 &+&\frac{1}{2}\sum_{i=1}^{d}\Gamma_{x_{s}}(O_{s}e_{i}).
\Gamma_{x_{s}}(O_{s}e_{i})
-\Gamma_{x_{s}}'(O_{s}e_{i}, O_{s}e_{i})\}. X_{s},\label{ito}
\end{eqnarray}
so that $(ii)$ follows by the boundedness of $X_{s}, \Gamma, \Gamma'$.

%%%%%%%%%%%%%%%%%%%%%%%%%%%%%%%%%%%%%%%%%%%%%%%%%%%%%%%%%%%%%%%%%%%%%%%%%%%

\section{Existence and uniqueness results for the transformation}
First, let us introduce some notations.

\noindent
The norms on vector spaces like $\R^{N},L(\R^{N}, \R^{m})$ will be denoted
by $|.|$. For a random variable $F$ we set $\|F\|_{\infty}=
\|F\|_{L^{\infty}(\Omega)}$.
 We denote by $BS(\R^{n})$ 
the set of Brownian semimartingales (terminology used in 
 \cite{D}) i.e.
the semimartingales $x_{s}$ such that $dx_{s}= O_{s}dw_{s}+A_{s}ds$ where 
$(w_{s})$ is an $\R^{m}$-valued Brownian motion, $(O_{s})$ is a $L(\R^{m},
\R^{n})$-valued adapted process and $(A_{s})$ is an $\R^{n}$-valued adapted
process. Similarly $BS(M)$ is the subset of $x_{.}\in BS(\R^{N})$ such that
$x_{s}\in M$ a.s. In the sequel we denote by $E$ the set of 
$x_{.}=x_{0}+\int_{0}^{.}O_{s}dw_{s}+\int_{0}^{.}A_{s}ds \in BS(\R^{N})$ 
such that
\[ \|x\|_{E}:=(E\sup_{s\in[0,1]}|O_{s}|^{2}+E\int_{0}^{1}|A_{s}|^{2}ds)^{1/2}
< +\infty.\]
$(E,\|.\|_{E})$ is a Banach space.

\subsection{Statement of the result}
\begin{defi}
Let $\sigma^{0}$ be a semimartingale on $M$. A map $\sigma: \R\longrightarrow
BS(M)$ is a solution of the initial value problem
\[
{\cal P}_{1}:
\left\{ \begin{array}{c}
\displaystyle \frac{d\sigma}{dt}=H(\sigma(t)) h(w)\\
\sigma(0)=\sigma^{0}.
\end{array}\right. 
\]
if $\sigma$, viewed as a map $\R\longrightarrow E$, satisfies
\[
{\cal P}_{2}:\left\{\begin{array}{c}
\displaystyle \frac{d\sigma}{dt}=\tH(\sigma(t)) h(w)\\
\sigma(0)=\sigma^{0}
\end{array}\right.
\]
($E$ is equiped with the norm $\|.\|_{E}$) and $\sigma_{s}(t)\in M, \; \; s\in [0,1],\; a.s.$.
\end{defi}
\begin{Rq}{\rm
In this definition we restrict ourselves to the solutions $\sigma$ such that
$\sigma(t)\in BS(M)$ for each $t$.}
\end{Rq}
\noindent
Let $x$ be a semimartingale in $E$ with $dx_{s}=O_{s}dw_{s}+A_{s}ds$. We will
say that $x$ satisfies the hypothesis (H) if

\[\|\sup_{s}|O_{s}|\|_{\infty}+\|\int_{0}^{1}|A_{s}|^{2}ds\|_{\infty}<\infty\]
We can now state the existence and uniqueness result concerning the
problem ${\cal P}_{1}$.
\begin{thm}
let $h:[0,1]\times W\longrightarrow \R^{d}$ be such that $h\in \H, a.s.$ and 
$\int_{0}^{1}\h_{t}^{2}(w)\leq C$ a.s. for some constant $C$. Then for every 
semimartingale $\sigma^{0}\in BS(M)\cap E$ which satisfies the hypothesis 
{\rm (H)} there exists a unique solution
$t\mapsto \sigma(t)$ defined on $\R$ to the problem ${\cal P}_{1}$.
\end{thm}
Let us point out that the main problem here consits in solving an ordinary
differential equation which does not satisfy the local Lipshitz condition.

\subsection{Proof of the theorem}
In this section $\sigma^{0}\in BS(M)$ and $h\in \H$  are fixed and satisfy 
the assumptions of Theorem 3.1. We will set $\|h\|_{\infty}:=
 \|\sup_{s}|h_{s}|\|_{\infty}$ (which is $\leq 
\|(\int_{0}^{1}|\h_{s}|^{2}ds)^{1/2}\|_{\infty} \leq C$). 

\noindent
Let $(x_{s})\in E$. The following inequality, which is a consequence of
the Burkholder inequalities, will be frequently used :
\begin{equation}
\label{bdg}
E\sup_{s\in[0,1]}|x_{s}|^{2} \leq C \|x\|_{E}^{2}.
\end{equation}
\noindent
Let ${\cal E}$ be the space of paths on $E$ i.e. the set of continuous maps
 $\sigma : \R\longrightarrow E$. We define the map 
\begin{eqnarray*}
\l : && \E\longrightarrow \E\\
&&  \sigma \mapsto \sigma_{0}+\int_{0}^{.}\tH(\sigma(t))h dt,
\end{eqnarray*}
where the last integral is a Riemann integral in the Banach space $E$.

\begin{lm}
let $x^{1}, x^{2}$ be in $E$, with $x^{i}_{.}=x^{i}_{0}+
\int_{0}^{.} O_{s}^{i}dw_{s}+\int_{0}^{.}A_{s}^{i}ds$, then
\begin{equation}
\label{lm1}
 E\sup_{s\in [0,1]}  | \tH_{s}(x^{1})- \tH_{s}(x^{2})|^{2}\leq
K_{1}(x^{1},x^{2})\|x^{1}-x^{2}\|_{E}^{2},
\end{equation}
where $K_{1}(x^{1},x^{2})$ is of the a continuous function of $\|\sup_{s\in [0,1]} 
|O_{s}^{i}|\|_{\infty},\| \int_{0}^{1}|A_{s}^{i}|^{2}\|_{\infty}, i=1,2$.
\end{lm}

\noindent
{\bf Proof.} We shall use the following notation:
 $X_{s}^{i}:= \tH_{s}(x^{i}),\; i=1,2,\; \Delta_{s}:= X_{s}^{1}-X_{s}^{2},\; 
\psi(t)=\sup_{\tau\leq t}|\Delta_{\tau}|^{2}.$ From the 
Proposition 2.1, $X_{s}^{i}, i=1,2$ satisfy equations of the form
\[ dX_{s}^{i}= G_{s}^{i} (dw_{s})X_{s}^{i} + F_{s}^{i} X_{s}^{i} ds,\; i=1,2. \]
 Hence 
\[ d\Delta_{s}= (G_{s}^{1} (dw_{s}) \Delta_{s} + (G_{s}^{1}-G_{s}^{2}) (dw_{s})X_{s}^{2})
+ (F_{s}^{1}\Delta_{s} +(F_{s}^{1}-F_{s}^{2})X_{s}^{2}) ds,\]
and 
\begin{eqnarray*}
 E\psi(t)&\leq& 4E\sup_{\tau\in[0,t]}\int_{0}^{\tau}[|G_{s}^{1}|^{2} 
|\Delta_{s}|^{2}ds + |X_{s}^{2}|^{2}|G_{s}^{1}-G_{s}^{2}|^{2}]ds \\
               &+& 4E\sup_{\tau\in[0,1]}|\int_{0}^{\tau}F_{s}^{1} 
                   \Delta_{s}ds|^{2} +4E\sup_{\tau\in [0,t]}\int_{0}^{\tau}
                      |X_{s}^{2}|^{2}|F_{s}^{1}-F_{s}^{2}|^{2}ds .
\end{eqnarray*}
But
\begin{eqnarray*}
E\sup_{\tau\in [0,t]}|\int_{0}^{\tau}F_{s}^{1} \Delta_{s}|^{2}ds &\leq& 
E\sup_{\tau\in [0,t]}[\int_{0}^{\tau}|F_{s}^{1}|^{2}ds 
\int_{0}^{\tau}|\Delta_{s}^{2}|ds]\\
                      &\leq& \|\int_{0}^{1}|F_{s}^{1}|^{2}ds\|_{\infty}
                                    \int_{0}^{t}E|\Delta_{s}|^{2}ds\\
  &\leq&\|\int_{0}^{1}|F_{s}^{1}|^{2}ds\|_{\infty}\int_{0}^{t}\psi(s)ds,
\end{eqnarray*}
which implies that
\begin{eqnarray*}
E\int_{0}^{t}|G_{s}^{1}|^{2} |\Delta_{s}|^{2}ds &\leq& \|\sup_{s\in[0,1]}
|G_{s}^{1} |^{2}\|_{\infty} \int_{0}^{t}E |\Delta_{s}|^{2}ds \\
 &\leq&\|\sup_{s\in[0,1]}|G_{s}^{1}|^{2}\|_{\infty}\int_{0}^{t}\psi(s)ds.
\end{eqnarray*}
Therefore
\begin{eqnarray*}
\psi(t)&\leq& 4( \|\int_{0}^{1}|F_{s}^{1}|^{2}ds\|_{\infty}+
\|\sup_{s\in[0,1]}|G_{s}^{1}|^{2}\|_{\infty})\int_{0}^{t}\psi(s)ds \\
               &+&4 \|X_{s}^{2}\|_{\infty}^{2}E\int_{0}^{1}[|G_{s}^{1}-
                G_{s}^{2}|^{2}+|F_{s}^{1}-F_{s}^{2}|^{2}]ds.
\end{eqnarray*}
On the other hand, using (\ref{ito}) and the boundedness of $\Gamma'$ we 
see that there is a constant $C_{1}$ such that
\begin{eqnarray*}
  E\int_{0}^{1}|G_{s}^{1}-G_{s}^{2}|^{2}ds &\leq &E\sup_{s}|G_{s}^{1}-
G_{s}^{2}|^{2} \\
&\leq & C_{1}\|\sup_{s}|O_{s}^{1}|^{2}\|_{\infty}(E\sup_{s}|O_{s}^{1}-
O_{s}^{2}|^{2}+ E\sup_{s}|x^{1}-x^{2}|^{2})\\
&\leq& C_{2}\|\sup_{s}|O_{s}^{1}|^{2}\|_{\infty} \|x^{1}_{s}-x^{2}_{s}
\|_{E}^{2}.
\end{eqnarray*}
Similarly, we have 
\[     E\int_{0}^{1}|F_{s}^{1}-F_{s}^{2}|^{2}ds\leq C_{3}(\|\sup_{s}|
  O^{i}_{s}|\|_{\infty}, i=1,2)\|x^{1}-x^{2}\|_{E}^{2},\]     
where $C_{3}$ is a polynomial function of $\|\sup_{s}| O^{i}_{s}|\|_{\infty}, 
i=1,2 $. Therefore
\[ \psi(t)\leq C_{4}(x^{1}, x^{2})\int_{0}^{t}\psi(s)ds +C_{5}(x^{1},x^{2})
\|x^{1}-x^{2}\|^{2}, \]
where $C_{4}, C_{5}$ are (polynomial) functions of $\|\sup_{s}| O^{i}_{s}|\|_{
\infty}, \|\int_{0}^{1}|F_{s}^{i}|^{2}, i=1,2 $. Now (\ref{lm1}) follows from 
the last inequality and the Gronwall lemma. 

\begin{lm}
For $x^{1}, x^{2}\in E$ such that $x^{1}_{s}, x_{s}^{2}\in Y$ a.s., we have 
\begin{equation}
\label{lm21}
\| \tH (x^{1})h-\tH(x^{2})h\|_{E}\leq K_{2}(x^{1}, x^{2})\|x^{1}-x^{2}\|_{E},
\end{equation}
where $K_{2}(x^{1},x^{2})$ is a continuous function of $\|\sup_{s}|O_{s}^{i}|\|_{\infty}, \|\int_{0}^{1}|A_{s}^{i}|^{2}ds \|_{\infty}, i=1,2.$ Consequently 
for $\sigma^{1}, \sigma^{2}\in \E$ :
\begin{equation}
\label{lm22}
 \| \l(\sigma^{1})(t)-\l(\sigma^{2})(t)\|_{E}\leq \sup_{s\in [0,t]}
K_{2}(\sigma^{1}(s), \sigma^{2}(s))\int_{0}^{t}\|\sigma^{1}(t)-\sigma^{2}(t)\|_{E}.
\end{equation}
\end{lm}
{\bf Proof.} As usual we set $X_{s}^{i}=\tH_{s}(x^{i}), i=1,2$ and
$dx_{s}^{i}=O_{s}^{i}dw_{s}+A_{s}^{i}ds, i=1,2$. We have
\begin{eqnarray*}
d(Xh)_{s} &=& -\Gamma_{x_{s}}(O_{s}dw_{s}).X_{s}h_{s}-[\Gamma_{x_{s}}(A_{s}).
X_{s}h_{s}\\
&+&\frac{1}{2}\sum_{i=1}^{N}\Gamma_{x_{s}}(O_{s}e_{i}).\Gamma_{x_{s}}(
O_{s}e_{i}).X_{s}h_{s} \\
&-& \frac{1}{2}\sum_{i=1}^{N}\Gamma_{x_{s}}'.(O_{s}e_{i}, O_{s} e_{i}).X_{s}
h_{s} +X_{s}\h_{s}]ds.
\end{eqnarray*}
Then, if we set $d\Delta_{s}:=o_{s}dw_{s}+a_{s}ds:=d(X^{1}h-X^{2}h)_{s}$
we have:
\begin{eqnarray*}
E\sup_{s\in [0,1]}|o_{s}|^{2}&\leq&\|h\|_{\infty}^{2}\{|\Gamma_{x_{s}^{1}}|^{2}
\|\sup_{s}|X_{s}|^{2}\|_{\infty} E\sup_{s}|O_{s}^{1}-O^{2}_{s}|^{2} \\
&+& \|\sup_{s}|O_{s}^{2}|^{2}\|_{\infty} E\sup_{s}|G(x_{s}^{1}, X_{s}^{1})-
G(x_{s}^{2}, X_{s}^{2})|^{2}\}
\end{eqnarray*}
where $G$ is a Lipshitz function (since $x_{s}\in Y$ a.s. and $X_{s}, 
\Gamma, \Gamma'$ are bounded). Hence 
\[  E\sup_{s}|G(x_{s}^{1}, X_{s}^{1})-G(x_{s}^{2}, X_{s}^{2})|^{2}\leq
C_{1} E \sup_{s}(\|x_{s}^{1}-x_{s}^{2}\|^{2}+\|X_{s}^{1}-X_{s}^{2}\|^{2}).
\]
Therefore, using lemma 3.1 and the Burkholder inequality (see (\ref{bdg})) we
 find that
\[  E\sup_{s\in [0,1]}|o_{s}|^{2} \leq C_{2}\|h\|_{\infty} (\|\sup_{s}
|O_{s}^{2}|^{2}\|_{\infty} +1)\|x^{1}-x^{2}\|_{E}^{2}. \]
We turn now to control the term $E\int_{0}^{1}|a_{s}|^{2}ds $. We have
\[ E\int_{0}^{1}|a_{s}|^{2}ds\leq \Delta_{1}+\Delta_{2}+\Delta_{3}+\Delta_{4},
\]
with
\begin{eqnarray*}
 \Delta_{1} &=& E\int_{0}^{1}|(\Gamma_{x_{s}^{1}}(A_{s}^{1})X_{s}^{1}-
\Gamma_{x_{s}^{2}}(A_{s}^{2})X_{s}^{2})h_{s}|^{2}ds\\
   &\leq& C_{3}\|h\|_{\infty}\{\int_{0}^{1}|A_{s}^{1}-A_{s}^{2}|^{2}ds\\ 
  &+& E\int_{0}^{1}|G(x_{s}^{1}, X_{s}^{1})-G(x_{s}^{2}, X_{s}^{2})|^{2}
|A_{s}^{2}|^{2} ds\\
 &\leq & C_{3}^{'} \|h\|_{\infty}\{ E\int_{0}^{1}|A_{s}^{1}-A_{s}^{2}|^{2}ds\\
 &+& \|\int_{0}^{1}|A_{s}^{2}|^{2}ds \|_{\infty}E\sup_{s}|G(x_{s}^{1},
 X_{s}^{1})-G(x_{s}^{2}, X_{s}^{2})|^{2}\},
\end{eqnarray*}
where $G$ is (another) Lipshitz function. By the same arguments as above
we get 
\[ \Delta_{1}\leq C_{3}^{''}\|h\|_{\infty}(\|\int_{0}^{1}|A_{s}^{2}|^{2}ds 
\|_{\infty}+1) \|x^{1}-x^{2}\|^{2}_{E}. \]
By similar majorizations, we get easily
\begin{eqnarray*}
\Delta_{2}&:=& E\int_{0}^{1}\frac{1}{2}|\sum_{i}\Gamma_{x_{s}^{1}}(O_{s}^{1}
e_{i}).\Gamma_{x_{s}^{1}}(O_{s}^{1}e_{i}). X_{s}^{1}h_{s}-\\
&-& \Gamma_{x_{s}^{2}}(O_{s}^{2}e_{i}).\Gamma_{x_{s}^{2}}(O_{s}^{2}e_{i}). 
X_{s}^{1}h_{s} |^{2}ds \\
&\leq& \|h\|_{\infty}C_{4}(\|\sup_{s}|O_{s}^{i}|\|_{\infty},i=1,2)\{ 
 E\sup_{s}|O_{s}^{1}-O_{s}^{2}|^{2}+\|x^{1}-x^{2}\|_{E}^{2} \}
\end{eqnarray*}
and
\begin{eqnarray*}
\Delta_{3}&:=&\int_{0}^{1}|\frac{1}{2}|\sum_{i}\Gamma_{x_{s}^{1}}'.(O_{s}^{1}e_{i}, O_{s}^{1}.e_{i}).X_{s}^{1}h_{s}^{1}-\\
&-& \Gamma_{x_{s}^{1}}'.(O_{s}^{1}e_{i}, 
O_{s}^{1}.e_{i}).X_{s}^{1}h_{s}^{1}|^{2}ds \\
&\leq & \|h\|_{\infty}C_{5}(\|\sup_{s}|O_{s}^{i}|^{2}\|^{\infty},i=1,2)
\{ E\sup_{s}|O_{s}^{1}-O_{s}^{2}|^{2}+\|x^{1}-x^{2}\|_{E}^{2} \}.
\end{eqnarray*}
where $C_{4}, C_{5}$ are (polynomial) functions of their arguments.
For the last term, we have by using Lemma 3.1
\begin{eqnarray*}
\Delta_{4}&:=& E\int_{0}^{1}|X_{s}^{1}-X_{s}^{2}|^{2}|\h_{s}|^{2}ds\\
    &\leq& \|\int_{0}^{1}|\h_{s}|^{2}ds\|_{\infty}E\sup_{s}
 |X_{s}^{1}-X_{s}^{2}|^{2}\\
  &\leq & K_{1}(x^{1}, x^{2})\|\int_{0}^{1}|\h_{s}|^{2}ds\|_{\infty}
\|x^{1}-x^{2}\|_{E}^{2}.
\end{eqnarray*}
The inequality (\ref{lm21}) is now clear and the second inequality of the 
lemma is immediate. $\Box$

\noindent
Now, given $T>0$ we denote by $\E_{T}$ de the space of maps 
$\sigma :[-T,T]\longrightarrow E$ such that if $d\sigma_{s}(t)= O_{s}(t)dw_{s}
+A_{s}(t)ds $ we have \\

\noindent
{\bf (i)} $\sup_{s\in [0,1]}|O_{s}(t)|^{2}+\int_{0}^{1}|A_{s}(t)|^{2}ds\leq 1 +
\|\sup_{s}|O_{s}^{0}|^{2}\|_{\infty}+\|\int_{0}^{1}|A_{s}^{0}|^{2}ds
\|_{\infty}$  a.s.,
\\
{\bf (ii)}$\sigma_{s}(t)\in Y $  a.s.

\begin{lm} 
For $T$ sufficiently small ($T\leq T_{0}$, say), we have $\l(\E_{T})\subset
 \E_{T}$. Furthermore $T_{0}$ depends only on $\|\sup_{s}|O_{s}^{0}|\|_{
\infty}, \|\int_{0}^{1}|A_{s}^{0}|^{2}ds\|_{\infty}$ and $\|\int_{0}^{1}|
\h_{s}|^{2}ds\|_{\infty}$. Finally, there is a unique map 
$\sigma: [-T_{0},T_{0}]\longrightarrow E$ which 
satisfies
\[ \sigma(t)=\sigma(0)+\int_{0}^{t}H(\sigma(s))h ds, \; \mbox{for all} \;
t\in [-T_{0},T_{0}]. \]
\end{lm}

\noindent
{\bf Proof.} Let $\sigma_{.}(t)\in \E_{T}$ with $d\sigma_{s}(t)=O_{s}(t)dw_{s}+
A_{s}(t)ds$  and $X_{s}(t)=\tH_{s}(\sigma(t)),\;z_{.}(t)=X_{.}(t)h_{.},\; Z(t)=\l(\sigma)(t)$. From the equation (\ref{ito}), we deduce easily that if we set
$dz_{s}(t)= o_{s}(t)dw_{s}+b_{s}(t)ds $ then 
\[ \sup_{s\in[0,1]}|o_{s}(t)|\leq C_{1} \sup_{s\in[0,1]}|O_{s}(t)| \]
\[ \int_{0}^{1}|a_{s}(t)|^{2}ds\leq C_{2}(\sup_{s\in[0,1]}(|O_{s}(t)|^{2}+
|O_{s}(t)|^{4})+
\int_{0}^{1}|A_{s}(t)|^{2}ds +\int_{0}^{1}|\h_{s}|^{2}ds)\]
for some constants $C_{1}, C_{2}$. If we denote $dZ_{s}(t)=
\bar{O}_{s}(t)dw_{s}+\bar{A}_{s}(t)ds$, then using the fact that 
 $\sigma\in \E_{T}$, we have 
\[ \sup_{s\in [0,1]}|\bar{O}_{s}(t)|\leq \int_{0}^{t}C_{1}dt+\|\sup_{s}
|O_{s}^{0}| \|_{\infty},\]
\[ \int_{0}^{1}|\bar{A}_{s}(t)|^{2}ds\leq C_{2}\int_{0}^{t}(2+\int_{0}^{1}
    |\h_{s}|^{2}ds)dt+\int_{0}^{1}|A_{s}^{0}(t)|^{2}ds ,\]
and we see that, for $T$ sufficently small, the condition (i) is 
satisfied.

\noindent
For condition (ii), we have 
\[ \| \sup_{s\in [0,1]}|\l(\sigma(t))-\sigma_{0}|\|_{\infty}\leq
K|h|_{\infty}|t|,\]
 by the boundedness of $\tH$. Since $\sigma_{s}(0)\in M$ for
$s\in [0,1]$ we see again that for $T$ sufficiently small we have
: a.s. $\l(\sigma(t))_{s}\in Y$ for $s\in [0,1]$. The first assertion of the 
lemma is now proved.

\noindent
The proof of the second assertion is standard in view of the first: we have 
a constant $K$ such that for all $\sigma^{1}, \sigma^{2}\in \E_{T_{0}}$:
\[ \|\l(\sigma^{1}(t))-\l(\sigma^{2}(t))\|_{E}\leq K\int_{0}^{t}\|\sigma^{1}(t)
-\sigma^{1}(t)\|_{E}dt .\]
The constant $K$ is now independent of $\sigma^{1}, \sigma^{2}$. By usual 
arguments one can verify the second assertion of the lemma.
Namely, for some $n_{0}$ sufficiently large, the map $\l^{(n_{0})}:\E_{T}
\longrightarrow \E_{T}$ is a contraction ($\E_{T}$ is endowed with the
norm \\
$\sup_{t\in[-T_{0},T_{0}]}\|.\|_{E}$). Hence $\l^{n_{0}}$ and then $\l$ admit
a unique fixed point $\sigma$. $\Box$

The next lemma garantees the boundedness of the solution to ${\cal P}_{1}$.
Its proof is a slight modification of the proof of [\cite{D}, Proposition 7.1].

\begin{lm}
Suppose that $t\mapsto\sigma(t)$ is a solution to ${\cal P}_{1}$, then there
exists a function $\beta$ such that for all $t\in \R$ :
\[ \|\sigma(t)\|_{E}\leq \beta(t)<+\infty.\]
\end{lm}

\noindent
{\bf Proof.} Let $\sigma_{.}(t)=\sigma_{0}(t)+\int_{0}^{.}O_{s}(t)dw_{s}+
\int_{0}^{.}A_{s}(t)ds $ be a solution to ${\cal P}_{1}$. By [\cite{D}, Lemma
7.2 ], $O(t),A(t)$ satisfy
\[
\left\{\begin{array}{c}
\displaystyle\frac{d O(t)}{dt}= C(\sigma(t))O(t)\\
\displaystyle\frac{d A(t)}{dt}= C(\sigma(t))A(t)+R(\sigma(t)).
\end{array}\right.
\]
where the derivatives are taken w.r.t the norms $\|O\|_{o}:=
E(\sup_{s\in[0,1]}|O_{s}|)$ and $\|A\|_{a}:=E(\int_{0}^{1}|A_{s}|^{2}ds)^{
1/2}$ and
$C$ is uniformely bounded (the bound, $K_{1}$ say, depends on $|h|_{\infty}$) 
and $R$ is of
the form $R(\sigma_{.}(t))=R^{1}(O_{.}(t))+H_{.}(\sigma(t))\h_{.}$, 
$R^{1}$ is a polynomial
function. Therefore, we have
for some constant $\gamma$ :
\[E |\frac{dO(t)}{dt}\|\leq \gamma\|O(t)\| \]
wich implies that 
\[E \|O(t)\|\leq \|O(0)\|e^{\gamma|t|}, \] 
here, as in \cite{D}, we use a ``vector'' version of 
the Gronwall lemma. Turning to $A(t)$, we have

\[E\int_{0}^{1}\left|\frac{d A_{s}(t)}{dt}\right|^{2}ds\leq 2K_{1}E\int_{0}^{1}
|A_{s}(t)|^{2}ds +4E\int_{0}^{1}|R^{1}_{s}(O(t))|^{2}ds +K_{2}E\int_{0}^{1}
\h_{s}^{2}ds, \]

hence for $\tau\leq t$
\[ \left\|\frac{dA(\tau)}{d\tau}\right\|\leq K_{3}(\|A(\tau)\|+e^{\gamma |t|}+ 
E(\int_{0}^{1}\h_{s}^{2}ds)^{1/2}), \]
and
\[ \|A(\tau)\|\leq e^{K_{3}|t|}[ \|A(0)\|+K_{3}|\tau|(
E(\int_{0}^{1}|\h_{s}|^{2})^{1/2}+e^{\gamma|t|})]. \]
We have used again the ``vector'' version of the Gronwall lemma.
$\Box$

\begin{lm}
Let $t\mapsto \sigma(t)$ be the solution to ${\cal P}_{2}$ in the interval
$[-T_{0}, T_{0}]$ (as constructed in lemma 3.3). Then there is a version of 
$\sigma$ such that a.s. the maps $(t,s)\mapsto \sigma_{s}(t)$ is continuously
differntiable in the $t$ variable.
\end{lm}
\noindent
{\bf Proof.} We will denote by $\dot{\sigma}(t)$ the derivative of $\sigma$
with respect to the norm $\|.\|_{E}$. For $t,t'\in [-T_{0},T_{0}]$ we have 
the estimates
\begin{eqnarray*}
E\sup_{s\in [0,1]} |\dot{\sigma}_{s}(t)-\dot{\sigma}_{s}(t')|^{2}&\leq& 
\||h|_{\infty}\|_{\infty}E\sup_{s\in [0,1]}
|H_{s}(\sigma(t))-H_{s}(\sigma(t'))|^{2}\\
        &\leq& \||h|_{\infty}\|_{\infty} K(\sigma(t))
                   \|\sigma(t)-\sigma(t')\|_{E}\\
         &\leq& \||h|_{\infty}\|_{\infty} K(\sigma(t))\int_{t}^{t'}\|
               H(\sigma(\tau)h\|_{E}d\tau\\
         &\leq& K'(\sigma(t)) |t-t'|^{2},
\end{eqnarray*}
where we have used lemma 3.1 in the second inequality. The last inequality 
is clear in view of the definition of $\E_{T_{0}}$; furthermore this 
definition 
implies that the constants $K(\sigma(t)), K'(\sigma(t))$ are independent of 
$t\in [-T_{0}, T_{0}]$, see lemma 3.1. We are now able to use a consequence
of the Kolmogorov lemma [\cite{D}, lemma 4.5.] which gives the desired result.
$\Box$

\noindent 
{\bf End of proof of the theorem:}

\noindent
$\bullet$ Since we have a local solution to the problem ${\cal P}_{2}$ let us
denote by $\tilde{\sigma}$ the maximal (unique) solution to this problem. By 
the boundedness lemma 3.4 this solution is defined on $\R$.

\noindent
$\bullet$ By lemma 3.5 we have a version of $\tilde{\sigma}$ such that 
$(t,s)\mapsto\tilde{\sigma}$ is continuously differentiable in the $t$ 
variable for $t\in [-T_{0}, T_{0}]$. Since $\R$ is a countable union of such
intervals we get a continuously differentiable version on the $t$ variable on
 $\R$. we set $\sigma_{s}(t)=\pi(\tilde{\sigma}_{s}(t))$ and using (i) of
Proposition 2.1 and the fact that $\tilde{\sigma}_{s}(0)\in M$ a.s., we see 
that $\sigma$ is a solution to the problem ${\cal P}_{1}$. 

\noindent
$\bullet$ By the uniqueness of the solution $\tilde{\sigma}$ to 
${\cal P}_{2}$ we have $\sigma=\tilde{\sigma}$ which shows the uniqueness of
the solution to ${\cal P}_{1}$ and completes the proof of the theorem. $\Box$

\section{Quasi-invariance of the Wiener measure on the path space}
Let us denote by $\Phi$ the development map, that is the map which 
associates
to an $M$-valued semimartingale $x_{s}=(\xi^{1}_{s},...,\xi^{d}_{s})$ 
the $\R^{d}$-valued semimartingale $\xi_{s}$ given by 
\[ \xi_{s}=\int_{X[0,s]}\theta\circ dX_{s}, \]
where $X_{s}$ is the horizontal lift of $x_{s}$, see Shigekawa \cite{Sh}.

In order to get the quasi-invariance property for the family of transformations $\sigma(t)$, it is necessary impose a condition to the torsion
$T$ of the connection. Namely, following \cite{D}, we say that the 
torsion  is 
``skew symetric'' or that the connection $\nabla$ is TSS if $g(T(X,Y),Y)=0$ 
for all vector fields $X,Y$, or, equivalently, if $v_{1}\mapsto \Theta_{u}(
v_{1},v_{2})$ is skew symetric for all $u\in O(M)$.

\begin{thm}
Let $\phi_{s}(t)=(\Phi(\sigma(t)))_{s}$ where $\sigma$ is the solution to $
{\cal P}_{1}$ with the initial condition $\sigma_{s}(0)=p_{s}$ with
$(p_{s})$ being the Brownian motion on $M$. Then $\phi(t)\in BS(\R^{d})$ and 
$t\mapsto \phi(t)$ is a 
solution to 
\begin{equation}
\label{eq4.1}
\frac{d\phi(t)(w)}{dt}=h_{.}(w)-\int_{0}^{.}\int_{0}^{s}
\Omega_{H_{\tau}(t)}(\circ d\phi_{\tau}(t), h_{\tau})\circ d\phi_{s}(t)
-\int_{0}^{.}\Theta_{H_{\tau}(t)}(\circ d\phi_{\tau}(t), h_{\tau}) 
\end{equation}
where the derivative is taken in $(BS(\R^{d}), \|.\|_{BS})$, $\|.\|_{BS}$ is 
the analogue of $\|.\|_{E}$ in $\R^{d}$ and $H_{.}(t)$
is the horizontal lift of $\sigma_{.}(t)$. Furthermore for each $t$, $\phi_{s}
(
t)$ is given by $d\phi_{s}(t)= o_{s}(t)dw_{s}+a_{s}(t)ds$ where $o_{s}(.)$ is
an adapted $O(d)$-valued process and $a_{s}(.)$ is an $\R^{d}$-valued 
process'such that for each $t$, there is a constant $C(t)$ s.t. :
$\int_{0}^{1}|a_{s}(t)
|^{2}ds \leq C(t)$ a.s. Therefore the law $\mu_{t}$ of the process
$s\mapsto \phi_{s}(t)$ is equivalent to $\mu$ and  the law $\nu_{t}$ is 
equivalent
 to $\nu$ with the same Radon-Nikod\'ym derivative : $d\nu_{t}/d\nu=
d\mu_{t}/d\mu$.  
\end{thm}

\noindent
{\bf Proof.} The proof of this theorem is an easy adaptation of the
 corresponding parts in \cite{D} or \cite{Hsu1}, taking into account  the fact that $h$ is in 
$\H$ and random and dealing with the appropriate norm.
First, using the structure equations (\ref{struct}), we prove that
$\phi(t)$ statisfies (\ref{eq4.1}). Next, write $d\phi_{s}(t)=o_{s}(t)dw_{s}
+a_{s}(t)ds $ (indeed, the development of a Brownian semimartingale 
is a Brownian semimartingale), then using (\ref{eq4.1}), one can show that
 $o(t), a(t)$ are solution to
\begin{eqnarray*}
 \frac{d o_{.}(t)}{dt}&= &c_{.}(\phi(t))o_{.}(t) \\
\frac{d a_{.}(t)}{dt} &=& c_{.}(\phi(t)) a_{.}(t) + \h_{.}\\
 &+& \frac{1}{2}({\rm ric}_{H_{.}(t)}h_{.}+\sum_{1}^{d}\Theta_{H_{.}(t)}'(
\e_{i},h_{.}, \e_{i}),
\end{eqnarray*} 
where the derivatives are taken w.r.t. the norms indicated in the theorem and $c_{s}(\phi(t))$ is the matrix given by
\[ c_{s}(\phi(t)) v= \int_{0}^{s}\Omega_{H_{\tau}(t)}(\circ d\phi_{\tau}(t),
 h_{\tau})v + \Theta_{H_{s}(t)}(h_{s}, v). \]
Hence, under the assumption on the torsion, $c_{s}(\phi(t))$ is skew-symetric
and $o_{s}(t)$ is orthogonal.  We omit the details and
refer to Theorem 5.1, Proposition 6.1 and Section 8 of \cite{D}. $\Box$

\begin{Rq} {\rm The anticipative case.

\noindent
As for the flat Wiener space, we can consider the case when $h$ is 
non-adapted (but, of course,  with other restrictions) and ask wether 
the problem ${\cal P}_{1}$ has a solution $\sigma(t)$ and if it has the 
quasi-invariance property. First, there are some difficulties in the definition
of the horizontal lift of an anticipative proceess (which requires to 
solve an anticipative SDE) and  the proof of the existence of the
solution to ${\cal P}_{1}$ has to be modified; for instance, we have
not an analogue to Burkholder inequalities. Second and most important, even 
if one has succeded to prove the existence of a solution $\sigma(t)$ to
${\cal P}_{1}$, we expect that the pullback $\xi(t)$ of $\sigma(t)$ to the 
flat Wiener space is of the form $d\xi_{s}(t)=o_{s}(t)\delta w_{s}+a_{s}(t)ds$,
 where $\delta$ is the Skorohod integral ; of course, one has to prove
an existence result for such equations. But the Wiener measure could hardly
be quasi-invariant under a transformation like $\xi(t)$. }
\end{Rq}

%%%%%%%%%%%%%%%%%%%%%%%%%%%%%%%%%%%%%%%%%%%%%%%%%%%%%%%%%%%%%%%%%%%%%%

\section{Anticipative integrals on a Riemannian manifold}
For the sake of simplicity, in this section $M$ is endowed with the 
Levi-Civita connection.
\subsection{Preliminaries}
{\it 5.1.1 }{\it Notations} : 
In this section $W=(C_{0}([0,1],\R^{d}), \H, \mu)$ will denote the Wiener 
space on $\R^{d}$, $D$ the usual Malliavin derivative on $W$ and $\delta$
its adjoint. We denote by ${\cal S}$ the set of smooth 
functionals $F$ on W, i.e. $F(w)=f(w_{s_{1}},....,w_{s_{n}})$ where
$f:(\R^{d})^{n}\longrightarrow \R$ is smooth. Recall that $p$ is the 
Brownian motion on $M$ determined 
by $p_{s}(w)=\pi(r_{s}(w))$ where $r_{s}$ satisfies $dr_{s}(w)=L_{i}(r_{s}(w))
\circ dw_{s}^{i}$, $r_{0}$ is given and the Ito stochastic parallel transport
is given by $t^{p}_{s_{1}\leftarrow s_{2}}= r_{s_{1}}\circ r_{s_{2}}^{-1}$.
The It\^o map is defined by $I:W\longrightarrow P(M)\; w\mapsto p(w)=\pi\circ 
r(w)$
. A tangent vector field on $P(M)$ is a process $\u(s)$ such $\u(s)\in T_{p(s)}
M$ for $s\in [0,1]$. 

\noindent
{\it 5.1.2 }{\it The gradient on the path space} (\cite{D}, \cite{Le}, 
\cite{FM}):

\noindent
$(a)$ Let $F:P(M)\longrightarrow \R$ be a 
cylindrical function i.e. $F(p)=f(p(s_{1}),..., p(s_{n}))$ where
$f:M^{n}\longrightarrow \R$ is a smooth function. Then, the gradient of
$F$ is the element $D^{M}F\in T(P)$ defined by 
\[ D_{s}^{M}F=\sum_{i=1}^{d}t^{p}_{s\leftarrow s_{i}}\nabla_{i}f 1_{s<s_{i}} ,
\]
where $\nabla_{i}f$ is the gradient w.r.t. the component $i$ of $f$ (defined
via the scalar product on $T_{p(s_{i})}M$). To each $h\in \H$ we associate
$h^{p}\in T(P)$ by $h^{p}(s)=t^{p}_{s\leftarrow 0}h(s)$ and we put
\[ D^{M}_{h}F= \int_{0}^{1}<D_{s}^{M}F, \h^{p}(s)>_{T_{p(s)}M}ds .\]

\noindent
$(b)$ This gradient is related to the transformation discussed in the above
paragraph as follows. For $h\in \H$ let $\sigma^{h}_{.}(t), t\in \R$ be 
the family of transformations defined in $\S 3$. Then for a cylindrical
function $F$ we have
\[ D^{M}_{h} F= \frac{d}{d\epsilon}|_{\epsilon=0} F(\sigma^{\epsilon h}(1)), \]
in $L^{2}(W)$. Indeed, it suffices to show that for $s\in [0,1]$, we have
\begin{equation}
\label{grad1}
 \frac{d}{d\epsilon}|_{\epsilon=0}\sigma^{\epsilon h}(1)= t^{p}_{s\leftarrow 
0}h(s).
\end{equation}
For this, let $q_{t}=(d/d\epsilon)_{\epsilon=0}(\sigma_{s}^{\epsilon h}(t))$.
Clearly, $q_{0}=0$ and using (\ref{eq1.1}) we get $(d/dt)q_{t}= H_{s}(p)h(s)$
which imlplies (\ref{grad1}).

\noindent
{\it 5.1.3 }{\it Integration by parts} (\cite{bi}, \cite{D}, \cite{FM}, 
\cite{Le}, \cite{AM}):

\noindent
 The integration by parts formula associated
to $D^{M}$ is the following (Bismut formula):
\begin{equation}
\label{ip1}
ED^{M}_{h}F= EF \int_{0}^{1}(\h(s)+\frac{1}{2}\ric h(s))dw_{s}.
\end{equation}
As a consequence of this formula, the operator $D^{M}$ is closable in $L^{2}(P)
\equiv L^{2}(W)$; we denote by $\D^{1,2}(M)$ its domain which is endowed with
the norm
\begin{equation}
\label{norm}
 \|F\|_{1,2}^{M}=\|F\|_{L^{2}}+ (E\int_{0}^{1}|D_{s}^{M}F|^{2}ds)^{1/2}.
\end{equation}
\noindent

{\it 5.1.4 }{\it The damped gradient} \cite{FM}:

\noindent
 The damped gradient, denoted by $\tD$, is
introduced in order to have an `ordinary' integration by parts formula 
instead of (\ref{ip1}) i.e.
\begin{equation}
\label{ip2}
E\tD_{h}F= EF\int_{0}^{1}\h(s)dw_{s}.
\end{equation}
Equivalently, $\tD$ must satisfy for every $F,h$ :$E\tD_{\h}F= ED^{M}_{\h_{1}}F$ where \\
$\h_{1}(s)+(1/2)\ric h_{1}(s) =\h$. After some calculations, one is led
to the following definition
\begin{equation}
\tD_{s} F= \sum_{i=1}^{n}r_{s}Q_{s,s_{i}}^{*}r_{s_{i}}^{-1}\nabla_{i}f 1_{s
<s_{i}} ,
\end{equation}
where $Q_{s,s'}:\R^{d}\longrightarrow \R^{d}$ is the solution of
\[ \frac{d Q_{s,s'}}{ds}=-\frac{1}{2}\ric Q_{s,s'},\; \; Q_{s',s'}=Id. \]
We then define an associated norm for $F$ by a formula like (\ref{norm}),
which, in turn, is equivalent to $\|.\|^{M}_{1,2}$, hence $D^{M}$ and
$\tD$ have the same domain.

\noindent
{\it 5.1.5 }{ \it Strong differentiability and the interwining formula} 
\cite{CM2}:

\noindent
A tangent process on $W$ is a process $\xi(s)$ which satisfies $\xi(0)=0$ and
$ d\xi(s)= A(s)dw_{s}+\h ds$ where an adapted $so(d)$-valued process and $\h$ is 
an adapted process in $L^{2}([0,1]\times W)$. For such a process and a 
cylindrical functional $F(w)=f(w_{s_{1}},....,w_{s_{n}})$ we define the 
derivative
\[ D_{\xi}F:= \sum_{i=1}^{n}\partial_{i}f.\xi(s_{i}) .\]
Observe that for such functionals $F\in {\cal S}$ we have
\begin{equation}
\label{sdiff1}
 D_{\xi}F=\frac{d}{d\epsilon}|_{\epsilon=0}F(\int_{0}^{.}e^{\epsilon A_{s}}
dw_{s} +\epsilon h_{s}) \; \mbox{ in} \; L^{2}(W) .
\end{equation}
Notation. For $d\eta(s)= A_{s}dw_{s}$ we set $D^{R}_{A}F:=D_{\eta}F$, 
 which means that the derivative corresponds to a rotation. Notice that
\[ D_{\xi}F= D_{h}F +D^{R}_{A}F. \]
The following integration by parts formula is an immediate consequence of 
(\ref{sdiff1}) and the invariance of the Wiener measure by rotation:
\[ ED_{\xi} F= EF\int_{0}^{1}\h(s)dw_{s}.\]
From that we deduce that $D_{\xi}$ is closable in $L^{2}(W)$. Denoting by 
$\D^{1,2}_{\xi}$its domain in $L^{2}(W)$, a random variable $F$ is said 
strongly differentiable iif $F\in\D^{1,2}_{\xi}$ for all $\xi$.

Similarly, given a cylindrical functional on $P(M)$ :$F(p)=f(p(s_{1}),...,
p(s_{n}))$ and $\xi$ a tangent process, we define
\[ D^{M}_{\xi}F:= \sum_{i=1}^{n}\partial_{i}f.\xi^{p}(s_{i}), \;
\mbox{with}\; \; \xi^{p}(s)=t^{p}_{s\leftarrow 0}\xi(s). \]
We define in the same way as above $\D^{1,2}_{\xi}(M)$ and the strong 
differentianility for $F:P(M)\rightarrow \R$. Then, for such a functional, 
we have : $F$ strongly diferrentiable iif $F\circ I$ is strongly 
differentiable and the following (interwining formula) holds
\begin{equation}
\label{intform}
D^{M}_{\xi}F= D_{\xi^{*}}(F\circ I),
\end{equation} 
with 
\[ d\xi^{*}(s)= d\xi(s)+\ro(s)\circ dw_{s},\; \; \ro(s)=\int_{0}^{s}\Omega_{
r(\theta)}(\circ dw_{\theta}, \xi(\theta)). \]

\subsection{ The anticipative integrals}
{\it Convention}. In this section, every process $u(s)\in T_{p(s)}M$ will
be identified to the process $t^{p}_{0\leftarrow s}u(s)$ also denoted by
$u(s)$. In particular, the gradients $D^{M}_{s}F, \tD_{s}F$ will designate
$t^{p}_{0\leftarrow s}D^{M}_{s}F$ and $t^{p}_{0\leftarrow s}\tD_{s}F$ 
respectively.
 Finally every random variable $F:P(M)\longrightarrow  \R$ will be identified 
with $F\circ I$.

\noindent
Notation. For a smooth functional $F$ and $i=1,...,d$ we define $\tD^{i} F$ by  

\[ \int_{0}^{1}\tD^{i}_{s}F. \h_{s}ds=
\int_{0}^{1}\tD_{s} F.\h_{s}
\e_{i} ds \] 
for every $h\in L^{2}([0,1])$ and if $F=(F_{1},...,F_{d})
\in (\tD^{1,2}(M))^{d}$ then we set :

\[ \tD.F:= \sum_{i=1}^{d}\tD^{i}F_{i}. \]
We define in the same way $D^{i} F$.

In  the following $\L^{1,p}$ will denote the set of $\R^{d}$-valued processes
$u$ such that $u(t)\in \D^{1,2}$ for almost all $t$ and
 
\[ \|u\|_{1,p}:= E(\int_{0}^{1}|u_{t}|^{2})^{p/2}+
E(\int_{0}^{1}\int_{0}^{1}|D_{s}u(t)|^{2}dsdt)^{p/2} < \infty. \]

Similarly, we denote by $\L^{1,p}(M)$ the set of the
processes $u(s)$ which satisfy the same condition when replacing $D$ by
$D^{M}$ (or $\tD$).  

\begin{defi} 
The anticipative integral associated to $D^{M}$ is the operator $\delta^{M}$
whose domain is
\[ {\rm Dom}(\delta^{M})=\{u\in L^{2}([0,1]\times W):\; (D^{M})^{*}u\in 
L^{2}(W)\}, \]
and defined by $\delta^{M}(u)= (D^{M})^{*}u$, where $(D^{M})^{*}$ is the
formal adjoint of $D^{M}$. Equivalently, we have
\[ {\rm Dom}(\delta^{M})=\{u\in L^{2}([0,1]\times W): \exists C(u), \forall
\phi\in \D^{2,1}_{M} : |E(\tD_{u}\phi)|\leq C(u)\|\phi\|_{2}\},\]
and for $u\in {\rm Dom}(\delta^{M})$, $\delta^{M}$ is the unique random 
variable in $L^{2}(W)$ which verifies
\[ E \int_{0}^{1}D_{s}^{M}\phi . u_{s}ds = E\phi\delta^{M}(u) ,\; {\rm for\;
 all}\; \phi\in \D^{2,1}_{M}. \]
Similarly, we define the anticipative integral $\deltam$ associated to the 
damped
gradient $\tD$ by replacing in the above formula $D^{M}$ by $\tD$.
\end{defi}

\noindent
The main result of Fang \cite{F}, is that if $p>2$ then $\L^{1,p}(M)\subset
Dom(\deltam)$.

\noindent
First we state the following proposition which is inspired by the Ogawa
method for defining noncausal integrals (see \cite{oga}, \cite{NZ}):

\begin{prop} 
Let $u\in \L^{1,2}(\R^{d})\cap \L^{1,2}(M)$ be a process such that a.s. 
the kernels 
$D_{s}u_{t}$ and $\tD_{s}u_{t}$ are of trace class. Then $u\in {\rm Dom}
(\deltam)$ and 

\begin{equation}
\deltam(u)= \delta(u)+ {\rm Trace}(D.u)-{\rm Trace}(\tD.u). 
\end{equation}
\end{prop}

\noindent
{\bf Proof.} Let $(\h_{i}), i\geq 1$ be an orthonormal basis of $L^{2}([0,1])$.
Without loss of generality, we suppose that $\u(s)=u(s)\e$ with $\e\in \R^{d}$. Then, we have $u=\sum_{1}^{\infty}\alpha_{i}(w)\h_{i}$ in $L^{2}
([0,1]\times W)$, with $\alpha_{i}=\int_{0}^{1}u_{s}\h_{i}(s)ds $.    
By the assumptions $\alpha_{i}\in \D^{1,2}\cap\D^{1,2}(M)$ and 
for $\phi \in \D^{2,1}(M)$, the integration by
part formula yields (we identify $\h_{i}$ with $\h_{i}\e$):

\begin{eqnarray*}
E\alpha_{i}\int_{0}^{1}\tD_{s}\phi .\h_{i}(s)ds&=& E\phi(
\alpha_{i}\int_{0}^{1}\h_{i}(s)dw_{s}- \tD_{h_{i}}
\alpha_{i})\\
 &=& E\phi(\int_{0}^{1}\alpha_{i}h_{i}(s)\delta w_{s}
+\int_{0}^{1}D_{s}\alpha_{i}.\h_{i}(s)ds -\int_{0}^{1}\tD_{s}
\alpha_{i}.\h_{i}(s)ds\\
&=& E\phi(\int_{0}^{1}\alpha_{i}\h_{i}(s)\delta w_{s}
+\int_{0}^{1}\int_{0}^{1}D_{s}u(t).\h_{i}(t).\h_{i}(s)ds\\ 
&-&\int_{0}^{1}
\int_{0}^{1}\tD_{s}u(t).\h_{i}(t).\h_{i}(s)ds ).
\end{eqnarray*}
Taking the sum on $i$, the first member converge to $E<\tD \phi, u>_{H}$,
and, under the assumptions of the proposition, the r.h.s. converges to
$\delta(u)+{\rm Trace}(D u)-{\rm Trace}(\tD u)$. $\Box$

\noindent
In \cite{K}, Kazumi follows the above method to calculate the 
Ornstein-Uhlenbeck operator associated to $\tD$ ($\tilde{L}:=-\deltam\circ
\tD$).

If $u(s)=\alpha(w) e(s)$ with $\alpha\in \D^{1,2}\cap
\D^{1,2}(M) $ and $e\in \H$ then $u\in {\rm Dom}(\deltam)$. For example, if
$u$ is a step process i.e.
\begin{equation}
\label{step}
 u(s)=\sum_{j=1}^{d}\sum_{j=1}^{n}\alpha_{i,j}(w)1_{[s_{i}, s_{i+1}[}
\e_{i},
\end{equation}
with $\alpha_{i,j}\in \D^{1,2}\cap\D^{1,2}(M)$ and $(s_{i}),i=1,...,n$ a
subdivision of $[0,1]$, then $u\in {\rm Dom}(\deltam)\cap\L^{1,2}\cap\L^{1,2}(M)$. For notational simplicity let $u(s)=\alpha(w)1_{[s_{1}, s_{2}[}\e$
with $\alpha\in \D^{1,2}\cap\D^{1,2}(M)$ and $\e \in \R^{d}$. Then, 
by the integration
by parts formula (\ref{ip2}), we have for $\phi\in {\cal S}$:
\begin{eqnarray*}
E\int_{0}^{1}\tD_{s}\phi. u(s)ds &=& E\phi(\int_{0}^{1}
\alpha 1_{[s_{1}, s_{2}[}(s)\e.\delta w_{s}+
\int_{0}^{1}D_{s}\alpha. \e 1_{[s_{1}, s_{2}[}(s)
ds \\
&-&\int_{0}^{1}\tD_{s}\alpha. \e 1_{[t_{1}, t_{2}[}(s)ds) \\
&=& \delta(u)+\int_{0}^{1}D_{s}.u_{s}ds -\int_{0}^{1}\tD_{s}.u_{s}ds.
\end{eqnarray*}
(In this case, the random variables $D_{s}.u(s), \tD_{s}.u(s)$ are clearly well
defined in $L^{2}(W)$ for $s\in [0,1]$). The above formula for $\deltam$
holds for every step process of the form (\ref{step}). The next theorem
extends this property to other processes.

\noindent
Following \cite{NP}, we say that a process $u$ belongs to $\L^{1,2}_{C}(M)$ if
 there is a version of $\tD u$ such that 

(i) The maps $s\mapsto \tD_{s\vee t}^{i}u_{s\wedge t}$ and  
$s\mapsto \tD_{s\wedge t}u_{s\vee t}$ are uniformly continuous (w.r.t $t$)
from $[0,1]$ onto $L^{2}(W,\R^{d}\times \R^{d})$.

(ii) $\sup_{s,t} E(|\tD_{s}u_{t}|^{2})<+\infty.$

Similarly, we define $\L^{1,2}_{C}$ by replacing $\tD$ by $D$ in the above
 definition. Then if $u\in \L^{1,2}_{C}(M)$, the processes 
\[ \tD^{+}_{t}.u_{t}:= \lim_{\epsilon \rightarrow 0, \epsilon >0} 
\sum_{i=1}^{d}\tD^{i}_{t}
u^{i}_{t+\epsilon}, \]
\[ \tD^{-}_{t}.u_{t}:= \lim_{\epsilon \rightarrow 0, \epsilon >0} 
\sum_{i=1}^{d}\tD^{i}_{t}
u^{i}_{t-\epsilon}, \]
exist in $L^{2}([0,1]\times W)$ and we can define similar processes 
for $u\in \L^{1,2}$

\begin{thm}
Let $u\in \L^{1,2}_{C}\cap\L^{1,2}_{C}(M)$. Then $u\in Dom(\deltam)$ and 
\begin{equation}
\label{formule}
\deltam(u)=\delta(u)+
\frac{1}{2}\int_{0}^{1}(D_{t}^{+}.u_{t}+
D_{t}^{-}.u_{t})dt-\frac{1}{2}\int_{0}^{1}(\tD_{t}^{+}.u_{t}+\tD_{t}^{-}.
u_{t})dt.
\end{equation}
\end{thm}

\noindent
{\bf Proof}. Let $u\in \L^{1,2}_{C}\cap\L^{1,2}_{C}(M)$. By definition
this implies that the r.h.s of (\ref{formule}) is an $L^{2}(W)$ random
variable; we will designate by $\deltam(u)$ this random variable. It remains
to prove that the following integration by parts formula holds
\begin{equation}
\label{eq5.1}
E\int_{0}^{1}\tD_{s}\phi.u(s)ds =E\phi\deltam(u) \; \;\mbox{for all}\;
 \phi\in {\cal S}.
\end{equation}
Let  
$(s_{i,n}=i/n), i=1,...,n$ be a subdivision of $[0,1]$. Define the step
process
\[ u^{n}(s)=\sum_{i=1}^{n}\frac{1}{n}\int_{s_{i,n}}^{s_{i+1,n}}u_{s}ds 
1_{[s_{i,n},s_{i+1,n}[}. \]
Then $u^{n}\in {\rm Dom}(\delta)\cap {\rm Dom}(\deltam)$ and we have
\[ \deltam(u^{n})=\delta(u^{n}) +A_{n} +\tilde{A_{n}}, \]
with
\[ A_{n}=\frac{1}{n}\sum_{i=1}^{n}\int_{s_{i,n}}^{s_{i+1,n}}\int_{s_{i,n}}^{
s_{i+1,n}}
D_{s}.u_{t}dsdt,\; \; \tilde{A}_{n}=\frac{1}{n}\sum_{i=1}^{n}
\int_{s_{i,n}}^{s_{i+1,n}}\int_{s_{i,n}}^{s_{i+1,n}}\tD_{s}.u_{t}dsdt. \]
This means that for every $\phi\in {\cal S}$ we have
\begin{equation}
\label{ipn}
 E\int_{0}^{1}\tD_{s}\phi.u_{s}^{n}ds = E\phi(\delta(u^{n})+A_{n}+
\tilde{A_{n}}) .
\end{equation}
Since $u\in \L^{1,2}$, we have  $u^{n}\longrightarrow u$ in $\L^{1,2}$ and 
$\delta(u^{n})\longrightarrow \delta(u)$ in 
$L^{2}(W)$ by Nualart-Pardoux \cite{NP}. Moreover we have
\[ E|A_{n}-\frac{1}{2}\int_{0}^{1}(D_{s}^{+}.u(s)+D_{s}^{-}.u(s))ds|
\longrightarrow 0,\]
see the proof of Theorem 3.1.1 p. 151 in Nualart \cite{N}. Similarly we
have
\[E|\tilde{A}_{n}-\frac{1}{2}\int_{0}^{1}(\tD_{s}^{+}.u(s)+\tD_{s}^{-}.u(s))ds|
\longrightarrow 0.\]
Then, using (\ref{ipn}) it follows that (\ref{eq5.1}) holds first for each 
bounded 
$\phi\in {\cal S}$ and then for each $\phi\in {\cal S}$. $\Box$

Let us denote by $\L^{k,p}$ the space $L^{p}([0,1],\D^{k,p})$. For the next
proposition we assume that $u$ satisfies the following: 
$u\in \L^{1,2}_{C}(M)$ and $u\in \L^{2,4}_{C}$ i.e. $u\in \L^{2,4}$ and

\[ s\mapsto D_{s\vee t}u_{s\wedge t}, s\mapsto D_{s\wedge t}u_{s\vee t} \]
 
are uniformly (in $t$) continuous from $[0,1]$ onto $L^{4}(W)$ and such that
$\sup_{s,t} E(|D_{s}u_{t}|^{4})< \infty$. 
We also assume that the processes $D^{+}_{t}. u_{t}+D^{-}_{t}.u_{t}$, 
$\tD^{+}_{t}. u_{t}+\tD^{-}_{t}.u_{t}$ are in $\L^{1,4}$.

\begin{prop} (It\^o formula) Let $u$ be a process which satisfies the 
above assumptions and 
 $\phi\in C^{2}(\R)$. Let $X_{t}:=\int_{0}^{t}u_{s}.\deltam p_{s}:=
\deltam(1_{[0,t]}u)$. Then
\[ \phi(X_{t})=\phi(X_{0})+ \int_{0}^{t}\phi'(X_{s})u_{s}.\deltam p_{s}-
\frac{1}{2}\int_{0}^{t}\phi''(X_{s})(\tD^{+}_{s}X_{s}+\tD^{-}_{s}X_{s}).u_{s}ds .\]
($\tD^{+}X_{s}$ is the vector $(\tD^{i,+}X_{s},i=1,...,d)$, with an obvious
 notation.)
\end{prop}
{\bf Proof.} Remark that
\begin{equation}
\label{ito2}
 X_{t}=\int_{0}^{t}u_{s}^{i}\circ dw_{s}^{i}-\frac{1}{2}\int_{0}^{1}
(\tD_{s}^{+}.u(s)+\tD_{s}^{-}.u(s))ds ,
\end{equation}
where the first integrals in the r.h.s are Stratonovich anticipative
integrals. The proposition follows easily from (\ref{ito2}) and the 
It\^o formula in the flat case, see Nualart [\cite{N},Theorem 3.2.3.].
$\Box$

\subsection{Remarks }
{\it 5.3.1.}  From 5.1.4. we see that a process $u$ belongs to 
${\rm Dom}(\delta^{M})$
iif $\tilde{u}\in {\rm Dom}(\deltam)$ where $\tilde{u}$ satisfy
\[ \tilde{u}_{s}+\frac{1}{2}\ric \int_{0}^{s}\tilde{u}_{\theta}d\theta=u(s),\]
and we have $\delta^{M}(u)=\deltam(u)$.

\noindent
{\it 5.3.2.}  Let $u(s)=\alpha(w)\h(s)$ where $\alpha$ is strongly 
differentiable and  $\h\in L^{2}([0,1],\R^{d})$. Then, by the 
integration by parts formula, it follows that $u\in {\rm Dom}(\delta^{M})$
and 
\begin{eqnarray*}
\delta^{M}(u)&=& \alpha\int_{0}^{1}(\h(s)+\frac{1}{2}\ric h(s))dw_{s}
-\int_{0}^{1}D^{M}_{s}\alpha. \h_{s}ds \\
 &=& \int_{0}^{1}\alpha(\h(s)+\frac{1}{2}\ric h(s))\delta w_{s}
 +\int_{0}^{1}D_{s}\alpha. (\h(s)+\frac{1}{2}\ric h(s))ds\\
&-& \int_{0}^{1}D^{M}_{s}\alpha. \h_{s}ds.
\end{eqnarray*}
But the interwining formula yields
\[ \int_{0}^{1}D^{M}_{s}\alpha. \h_{s}ds=\int_{0}^{1}D_{s}\alpha.
 (\h(s)+\frac{1}{2}\ric h(s))ds +D^{R}_{\ro_{h}}\alpha, \]
where $\ro_{h}(s)=\int_{0}^{s}\Omega_{r(s)}(\circ dw_{s}, h(s))$. Consequently
\begin{equation}
\delta^{M}(u)=\int_{0}^{1}(u(s)+\frac{1}{2}\ric (\int_{0}^{s}u(\theta)d\theta)) \delta w_{s} -D^{R}_{\ro_{h}}\alpha .
\end{equation}

\noindent
{\it 5.3.3.}  Now let $\u(s)$ be a process such that $\u\in \L^{1,2}$ and $\u(s)$ is
strongly differentiable for $s\in [0,1]$. We denote by $\h_{n}$ an orthonormal
basis of $L^{2}([0,1])$ and we assume without loss of generality
that $\u(s)=u(s)\e, \; \e\in \R^{d}$. Let $\u_{n}=\sum_{1}^{n}\alpha_{i}
\h_{i}\e$,
with $\alpha_{i}=<u,\h_{i}>_{L^{2}([0,1])}$.
We have $\u_{n}\longrightarrow u$ in $\L^{1,2}$ as $n\rightarrow\infty$. 
Also $\u_{n}\in {\rm Dom}(\delta^{M})$ and
\[ \delta^{M}(\u_{n})= \int_{0}^{1}(\u_{n}(s)+\frac{1}{2}\ric (\int_{0}^{s}
u_{n}(\theta).d\theta)) \delta w_{s}-\sum_{i=1}^{n}D^{R}_{\ro_{h_{n}}}\alpha.
\]
Then we can prove the following
\begin{prop}
(i) There exists a random variable in $L^{1}(W)$ which we denote also
by $\delta^{M}(\u)$ such that the integration by parts formula
$E<D^{M}\phi, u>= E\phi\delta^{M}(\u)$ holds for every smooth and bounded
functional $\phi$.

(ii)There is a constant $C>0$ such that
\begin{equation}
\label{estime-l1}
\|\delta^{M}(u)\|_{L^{1}}\leq \|u\|_{\L^{1,2}}+C(|\Omega|+|\Omega'|)
(\|u\|_{\L^{1,2}}+\|u\|_{\L^{1,2}}^{2})
\end{equation}
where $|\Omega|=\sup_{r}|\Omega_{r}|, |\Omega'|=\sup_{r}|\Omega_{r}'|  $.
\end{prop}
We omit the proof of these facts here, which is not difficult (we begin 
by the case where the $\alpha_{i}$ are cylindrical).

\noindent
{\it 5.3.4.}  We can do the same remark for the integral $\deltam$: As in 5.3.2, let
$\u(s)(w)=\alpha(w)\h_{s}$. then we have :

\[ \deltam(\u)=\delta(\u)+ <D\alpha, \h>-<\tD\alpha,\h>.\] 
But 
\[<\tD\alpha,\h>=<D^{M}\alpha,\tilde{\u}>  \]

with the notation of 5.3.1, then the interwining formula yields

\[ <\tD\alpha, \u>= < D\alpha, \u>- D^{R}_{\ro_{\tilde{h}}}\alpha \]
and
\[ \deltam(\u)=\delta(\u)+D^{R}_{\ro_{\tilde{h}}}\alpha .\]
and we have the same estimate as (\ref{estime-l1}). In view of the 
above expressions of $\deltam, \delta^{M}$, it seems not possible 
to obtain good estimates on the anticipative integrals by using the
Sobolev norms associated to the gradient in the flat Wiener space. This is 
essentially due to the rotational derivatives $D^{R}_{\ro}\alpha$
which appear in the expressions. For if $\gamma(s)$ is an adapted 
$so(d)$-valued process, then one can check easily that
\[ D^{R}_{\gamma}\alpha =\int_{0}^{1}D_{s}\alpha. \gamma(s)\delta w(s)+
 \int_{0}^{1} \gamma_{ij}(s) D_{s}^{i}D_{s}^{j}\alpha ds .\]
Hence the expressions of the anticipative integrals of a process $u$ involve
the second derivatives of $u(t)$.
We have similar difficulties with the norms associated to $\tD, D^{M}$.
 Let us finally mention that 
Cruzeiro and Fang (\cite{CF}, \cite{CF1}) have constructed a norm of the type $\|.\|_{1,2}$ in the 
space of tangent processes for which we have an estimate of the 
form $\|\delta^{M}(u)\|_{L^{2}}\leq c\|u\|_{1,2}$. Another attempt
in the case of Lie groups, to obtain an energy equality, has been made by 
Fang and Franchi \cite{FF}. 

\end{document}